\documentclass[12pt]{article}

\usepackage[left=1in,right=1in,top=1in,bottom=1in]{geometry}
\usepackage{mdframed}

\usepackage{subfigure}
\usepackage{color}

\usepackage{url,latexsym,amsmath,amssymb, amsfonts,enumerate,
  epsfig, graphics, amsthm}
  
\usepackage[colorlinks]{hyperref}

\newtheorem{theorem}{Theorem}
\theoremstyle{definition}
\newtheorem{example}[theorem]{Example}

\newtheorem{assumption}[theorem]{Assumption}

\newcommand\bbassump{\begin{mdframed}
 \begin{assumption}}
 
\newcommand\ebassump{\end{assumption}\end{mdframed}}

\newtheorem{algorithm}{Algorithm}

\newcommand{\EE}{\mathbb{E}}

\newcommand{\R}{{\mathbb R}}
\newcommand{\Z}{{\mathbb Z}}
\newcommand{\eqdef}{\overset{\text{def}}{=}}
\newcommand{\E}{\mathbb E}

\newcommand{\QNl}{\widehat{Q}^N}

\newcommand{\DN}{D^N}

\newcommand\hhf{{\scriptstyle{1\over 2 }\scriptstyle}}

\title{Computational complexity analysis for Monte Carlo approximations of classically scaled population processes}

\author{
David F. Anderson\thanks{Department of Mathematics, University of
  Wisconsin, Madison, USA.  anderson@math.wisc.edu, grant support from NSF-DMS-1318832 and Army Research Office W911NF-14-1-0401.},
\and
Desmond J. Higham\thanks{Department of Mathematics and
 Statistics, University of Strathclyde, UK.
d.j.higham@maths.strath.ac.uk, grant support from a Wolfson/Royal Society Research Merit 
Award  and an EPSRC Fellowship with reference 
EP/M00158X/1.
},
  \and
 Yu Sun\thanks{Facebook, Inc., 1 Hacker Way, Menlo Park, CA, USA.  ysunfb@fb.com.}
}

\begin{document}

\maketitle

\begin{abstract}
We analyze and compare the computational complexity of different simulation strategies for Monte Carlo in the setting of classically scaled population processes.  This allows a range of widely used competing strategies to be judged systematically.  Our setting includes stochastically modeled biochemical systems.  We consider the task of approximating the expected value of some path functional of the state of the system at a fixed time point.  We study the use of standard Monte Carlo when samples are produced by exact simulation and by approximation with tau-leaping or an Euler-Maruyama discretization of a diffusion approximation.  Appropriate modifications of recently proposed multilevel Monte Carlo algorithms are also studied for the tau-leaping and Euler-Maruyama approaches.  In order to quantify computational complexity in a tractable yet meaningful manner, we consider a parameterization that, in the mass action chemical kinetics setting, corresponds to the classical system size scaling.  We base the analysis on a novel asymptotic regime where the required accuracy is a function of the model scaling parameter. Our new analysis shows that, under the specific assumptions made in the manuscript,  if the bias inherent in the diffusion approximation is smaller than the required accuracy, then multilevel Monte Carlo for the diffusion approximation is  most efficient, besting multilevel Monte Carlo with tau-leaping by a factor of a logarithm of the scaling parameter.  However, if the bias of the diffusion model is greater than the error tolerance, or if the bias can not be bounded analytically, multilevel versions of tau-leaping are often the optimal choice.
\end{abstract}

\section{Introduction}
\label{sec:intro}

For some  large $N_0>0$ we consider a continuous time Markov chain satisfying the stochastic equation
\begin{equation}\label{eq:model_N0}
	X^{N_0}(t) = X^{N_0}(0) + \sum_{k=1}^K \frac1{N_0} Y_k\left(N_0\int_0^t \lambda_k(X^{N_0}(s)) ds\right) \zeta_k,
\end{equation}
where $X^{N_0}(t)\in \R^d$,  $K<\infty$, the $Y_k$ are independent unit Poisson processes and, for each $k$,  $\zeta_k \in \R^d$ and $\lambda_k:\R^d \to \R_{\ge 0}$ satisfies  mild regularity conditions.  
For a given path functional $f$, we consider the task of numerically approximating $\E[f(X^{N_0}(\cdot))]$, in the sense of confidence intervals, to some fixed tolerance $\varepsilon_0<1$.   Specifically, we consider the computational complexity, as quantified by the number of random variables utilized, required by different Monte Carlo schemes  to achieve a root mean squared error of $\varepsilon_0$.   For concreteness, we will assume throughout that the path functional $f$ depends upon $X^{N_0}(\cdot)$ only on the compact time interval $[0,T]$.

The class of models of the form (\ref{eq:model_N0}) 
has a long history in terms of modelling
 \cite{Bar58,Bar59,Del40,McQuarrie67},
analysis
\cite{AK2015,Kurtz72,Kurtz78} and computation
\cite{Gill76,Gill77}.
The framework covers many application areas, including population dynamics
\cite{Renshaw2011},
queueing theory \cite{Ross2006},
and several branches of
physics
\cite{Gardiner2002}. In recent years,
chemical and
biochemical kinetics
models in systems biology
\cite{Wilkinson2011}
have been
the driving force behind a resurgence of activity
in algorithmic developments, including
 tau-leaping
\cite{Gill2001} and its multilevel extension
\cite{AndHigh2012,AHS2014}.
In this setting, the
parameter $N_0$ in (\ref{eq:model_N0})
can represent Avogadro's number multiplied by the volume, and in this classical scaling, species are
measured in moles per liter.   More generally, however, $N_0$ can just be considered a large number, often of the order 100s or 1000s.

In section \ref{sec:summary}, we discuss some of the issues involved in quantifying computational complexity in the present setting, and introduce a novel scaling regime in which clear-cut comparisons can be made.   Further, the specific assumptions utilized throughout the manuscript are presented and a high-level summary of our main conclusions is presented.  
In section \ref{sec:approximations}, we summarize two widely used approximation methods for the 
 model \eqref{eq:model_N0}: the tau-leap discretization method, and the Langevin or diffusion approximation.  In section \ref{sec:CA}, we quantify the computational complexity of using exact simulation, tau-leaping, and simulation of the diffusion equation with standard Monte Carlo for approximating $\E[f(X^{N_0}(\cdot))]$ to a desired tolerance under our assumptions.  Further, in subsection \ref{sec:MLMC} we  review the more recent multilevel methods and quantify the benefits of their use in both the tau-leaping and diffusion scenarios.    In section \ref{sec:comp}, we provide numerical examples demonstrating our main conclusions.  In section \ref{sec:conc}, we close with some brief conclusions.

This paper makes use of results from two recent papers.
\begin{itemize}
\item  In \cite{AHS2014} an analysis was carried out to determine the variance of the difference between coupled paths in the jump process setting under a more general scaling than is considered here.
\item In \cite{AHS2016} an analysis was carried out to determine the variance of the difference between coupled paths in the setting of stochastic differential equations with small noise.  
\end{itemize}
Our goals here are distinct from those of these two papers.  First, the analysis in \cite{AHS2014} allowed such a general scaling that no modified versions of Euler based tau-leaping, such as midpoint or trapezoidal tau-leaping, could be considered.  Here, we consider a particular scaling (which is the most common in the literature) and present a unified computational complexity analysis for a range of Monte Carlo based methods.  This allows us to make what we believe are the first
concrete conclusions pertaining to the relative merits of current methods in a practically relevant asymptotic regime.  
Moreover, an open question in the literature involves the selection of the finest time-step in the unbiased version of multilevel Monte Carlo (since it is not constrained by the accuracy requirement).  By carrying out our analysis in this particular scaling regime, we are able to determine the asymptotics for the optimal selection of this parameter.  Selecting the finest time-step according to this procedure is shown to lower the computational complexity of the method by a nontrivial factor.  See the end of Section \ref{subsec:tauMLMC} for this derivation and the end of Section \ref{sec:comp} for a numerical example.  
Second, it has become part of the ``folk-wisdom'' surrounding these models that in the particular scaling considered here, properly implemented numerical methods applied to the diffusion approximation are the best choice.  This idea was somewhat exacerbated by the analysis in \cite{AHS2016}, which applied to a key aspect of the algorithm.   There it was shown that the variance between the coupled paths of a diffusion approximation is asymptotically smaller than the variance between the properly scaled jump processes.  However, here we show that the actual difference in overall complexity between properly implemented versions of multilevel Monte Carlo for the diffusion approximation and for the jump process never differ by more than a logarithm term.  If one combines this conclusion with the fact that the bias of the diffusion approximation itself is often unknown, whereas multilevel Monte Carlo applied to the jump process is naturally unbiased,  then the folk-wisdom is overturned and unbiased multilevel Monte Carlo is seen as a competitive choice.

\section{Scaling, assumptions, and a summary of results}
\label{sec:summary}

In order to motivate our analysis and
computations, we begin with a brief,
high-level, overview.
In particular, we discuss the
entries in Table~\ref{tab1}, which summarizes
the key conclusions of this work.
Full details  are given later in the manuscript, however we point out  here that the terms in Table \ref{tab1} include assumptions on the variances of the constituent processes that will be detailed below.

A natural approach to approximate the desired expectation is to
simulate paths exactly,
for example with the stochastic simulation algorithm
\cite{Gill76,Gill77}
or the next reaction method
\cite{Anderson2007a,Gibson2000},
in order to obtain independent sample paths
$\{ X^{N_0}_{[i]}  \}_{i=1}^{n}$
that can be combined into a sample average
\begin{equation}
 \hat \mu_n =  \frac{1}{n} \sum_{i=1}^{n} f(X^{N_0}_{[i]}(\cdot) ).
\label{eq:sampave}
\end{equation}
This becomes problematic if the cost of each sample path
is high---to follow a path exactly we must take account of
each individual transition in the process. This is  a serious issue
when many jumps take place, which is the case when $N_0$ is large.

The essence of the Euler tau-leaping approach is to fix
the system intensities over time intervals of length $h$,
and thereby only require the generation of $K$ Poisson random variables per time interval
\cite{Gill2001}.
In order to analyse the benefit of tau-leaping,
and related methods,
Anderson, Ganguly, and Kurtz
\cite{AndersonGangulyKurtz}
considered  a family of models, parameterized by $N\ge N_0$ (see \eqref{eq:model_family} below), and considered the limit
$N \to \infty$ and
$h \to 0$
with  $h= N^{-\beta} $ for some $\beta > 0$.
To see why such a limit is useful we
note two facts:
\begin{itemize}
\item If, instead, we allow $N \to \infty$ with $h$ fixed, then the stochastic fluctuations become negligible
\cite{AndKurtz2011,Kurtz72}.
In this \emph{thermodynamic limit} the model reduces to a deterministic
ODE, so a simple deterministic numerical method could be used.
\item If, instead, we allow $h \to 0$ with $N_0$ fixed then tau-leaping  becomes arbitrarily
    inefficient. The ``empty'' waiting times between
  reactions, which have nonzero expected
  values, are being needlessly refined by the discretization
  method.
\end{itemize}
The relation $h= N^{-\beta}$
brings together the large system
size effect (where exact simulation is expensive and
tau-leaping offers a computational
advantage) with  the small $h$ effect
(where the accuracy of tau-leaping can be analysed).
This gives a realistic setting where
the benefits of
tau-leaping can be quantified.
It may then be shown
\cite[Theorem~4.1]{AndersonGangulyKurtz}
that the bias arising from Euler tau-leaping is
$O(h) = O(N^{-\beta})$ in a wide variety of cases.  Higher order
alternatives to the original tau-leaping
method
\cite{Gill2001}
 are
available. For example, a
mid-point discretization
\cite[Theorem~4.2]{AndersonGangulyKurtz} or a trapezoidal method \cite{AndMaso2011} both achieve
$O(h^2) = O(N^{-2\beta})$ bias for a wide variety of cases.

As an alternative to tau-leap discretizations, we could replace the continuous-time Markov chain by a diffusion approximation
and
use
a numerical stochastic differential equation
(SDE) simulation method to generate
approximate paths \cite{AK2015}.  This approximation is detailed in section \ref{subsec:diffusion} below.  While higher order methods are available for the simulation of diffusion processes, we restrict ourselves to Euler-Maruyama as the perturbation in the underlying model has already created a difficult to quantify bias.  Thus, higher order numerical schemes are hard to justify in this setting.

For our purposes, rather than the step size $h$ of a particular approximate method, it is more natural to work in terms of the system size, $N_0$, and accuracy parameter $\varepsilon_0$.   Let $\varepsilon_0 = N_0^{-\alpha}$, for some fixed $\alpha > 0$.  A larger value of $\alpha$ corresponds to a more stringent accuracy requirement.
Next, consider the following family  of models parameterized by $N\ge N_0$,
\begin{equation}\label{eq:model_family}
	X^N(t) = X^N(0) + \sum_{k=1}^K \frac1N Y_k\left(N\int_0^t \lambda_k(X^N(s)) ds\right) \zeta_k,
\end{equation}
with initial conditions satisfying $\lim_{N\to \infty}  X^N(0) = x_0 \in \R^d_{>0}$.
We will study the asymptotic behavior, as $N\to \infty$, of the computational complexity required of various schemes to approximate   $\E[f(X^N(\cdot))]$ to a tolerance of 
\begin{equation}\label{eq:accsys}
\varepsilon_N= N^{-\alpha},
\end{equation}
where $f$ is a desired path functional.
Specifically, we require that both the bias and standard deviation of the resulting estimator is less than $\varepsilon_N$.

We emphasize at this stage that we are no longer studying
a fixed model. Instead we look at the
family of models
(\ref{eq:model_family})
parameterized through the system size $N$, and
consider the limit, as $N\to \infty$, of the computational complexity of the different methods
under the
accuracy requirement (\ref{eq:accsys}).
The computed results then tell us, to leading order, the costs associated with solving our fixed problem \eqref{eq:model_N0} with accuracy requirement $N_0^{-\alpha}$.

\subsection{Specific assumptions and a brief summary of results}

Instead of giving specific assumptions on the intensity functions $\lambda_k$ and the functional $f$, we give assumptions pertaining to the cost of different path simulation strategies, the bias of those strategies, and the variance of different relevant terms.  We then provide citations for when the assumptions are valid.  We expect the assumptions to be valid for a wider class of models and functionals than has been proven in the literature, and discovering such classes is an active area of research.

To quantify computational complexity, 
we define the ``expected cost-per-path'' to be the expected value of the number of random variables generated in the simulation of a single path.  
 Standard $\Theta$ notation is used (providing an asymptotic upper and lower bound in $N$ or $h$).  We emphasize that computations take place over a fixed time interval $[0,T]$. 
\begin{assumption}\label{assump:costperpath}
We assume the following  expected cost-per-path  for different methods.
\begin{table}[h]
\begin{tabular}{|c|c|}\hline
\textit{Method} & \textit{Expected cost-per-path}\\\hline
Exact simulation & $\Theta(N)$ \\\hline
Euler tau-leaping & $\Theta(h^{-1})$\\\hline
Midpoint tau-leaping & $\Theta(h^{-1/2})$\\\hline
Euler-Maruyama for diffusion  & $\Theta(h^{-1})$\\\hline
\end{tabular}
\end{table}
\end{assumption}

We make the following assumptions on the bias, $| \E[f(X^{N}(\cdot))] - \E[f(Z^N(\cdot))]|$, of the different approximation methods,  where $Z^N$ is a generic placeholder for the different methods.  

\begin{assumption}\label{assump:bias}
We assume  the following  biases.  
\begin{table}[h]
\begin{tabular}{|c|c|c|}\hline
\textit{Method} & \textit{Bias} & \textit{Reference}\\\hline
Exact simulation & 0 & N.A. \\\hline
Euler tau-leaping & $\Theta(h)$ & \cite{AndersonGangulyKurtz}\\\hline
Midpoint tau-leaping & $\Theta(h^{2})$ & \cite{AndersonGangulyKurtz} \\\hline
Euler-Maruyama for diffusion  & $\Theta(h)$ &\cite{AHS2016,KloedenPlaten92}\\\hline
\end{tabular}
\end{table}
\end{assumption}
A bias of $\Theta(h)$ for Euler-Maruyama applied to a diffusion approximation is extremely generous, as it assumes that the bias of the underlying diffusion approximation is negligible.
However, analytical results pertaining to the bias of the diffusion approximation for general functionals $f$ are sparse. 
 A startling result of the present analysis is that even with such generosity, the complexity of the \textit{unbiased} version of multilevel   tau-leaping is still often within a factor of a logarithm of the complexity of the multilevel version of Euler-Maruyama applied to the diffusion approximation.

We provide our  final assumption, pertaining to the  variances of relevant terms.  Below, $Z_h^N$ is a tau-leap process with step size $h$, ${\cal Z}_h^N$ is a midpoint tau-leap process with step size $h$, and $\DN_h$ is an Euler-Maruyama approximation of the diffusion approximation with step size $h$.  The coupling methods utilized are described later in the paper.  Finally, $h_\ell = M^{-\ell}$ for some integer $M>1$.
\begin{assumption}\label{assump:variance}
We assume the following relevant variances per realization/path. 

\vspace{.2in}

\noindent \begin{tabular}{|c|c|c|}\hline
\textit{Method} & \textit{Variance} &  \textit{Reference}\\\hline
Exact simulation & $\textsf{Var}(f(X^N(\cdot))) = \Theta(N^{-1})$& \cite{AHS2014} \\\hline
Euler tau-leaping & $\textsf{Var}(f(Z_h^N(\cdot))) = \Theta(N^{-1})$ & \cite{AHS2014} \\\hline
Coupled exact/tau-leap & $\textsf{Var}(f(X^N(\cdot)) - f(Z_h^N(\cdot))) = \Theta(h\cdot N^{-1})$ & \cite{AHS2014}\\\hline
Coupled tau-leap & $\textsf{Var}(f(Z_{h_\ell}^N(\cdot)) - f(Z_{h_{\ell-1}}^N(\cdot))) = \Theta(h_\ell\cdot N^{-1})$ & \cite{AHS2014}\\\hline
Midpt.~or trap.~tau-leaping & $\textsf{Var}(f({\cal Z}_h^N(\cdot))) = \Theta(N^{-1})$ & \cite{AHS2014} \\\hline
Euler-Maruyama for diffusion  & $\textsf{Var}(f(D_h^N(\cdot))) = \Theta(N^{-1})$ & \cite{AHS2016}\\\hline
Coupled diffusion approx. & $\textsf{Var}  ( f(\DN_{h_\ell}(\cdot)) - f(\DN_{h_\ell-1}(\cdot)))= \Theta(N^{-1}h_\ell^2+N^{-2}h_\ell)$ & \cite{AHS2016}\\\hline
\end{tabular}
\end{assumption}

The results presented in Table \ref{tab1} can now start coming into focus.  For example, we immediately see that in order to get both the bias and standard deviation under control, i.e.~below $\varepsilon_N$, we have the following:
\begin{description}
 \item[Monte Carlo plus exact simulation:] We require 
       $\Theta(N^{-1} \varepsilon_N^{-2}+1)$ paths for the standard deviation to be order $\varepsilon_N^2$, at
       a cost of $\Theta(N)$ per path.  This totals a computational complexity of $\Theta(\varepsilon_N^{-2}+N)$ or
          $\Theta(N^{2 \alpha}+N)$.
 \item[Monte Carlo plus tau-leaping:]
       $\Theta(N^{-1} \varepsilon_N^{-2}+1)$ paths at a cost of 
          $\Theta(\varepsilon_N^{-1})$ per path (required to achieve a bias of $O(\varepsilon)$), totaling a computational complexity of
          $\Theta(N^{-1}\varepsilon_N^{-3} + \varepsilon_N^{-1})$ or $\Theta(N^{3 \alpha-1}+N^{\alpha})$,
\end{description}
as summarized in the first two rows of Table~\ref{tab1}. Note  that the ``+1'' terms above
account for the requirement that we cannot generate less than one path.
In this regime, we see that tau-leaping is beneficial for
$\alpha < 1$. This makes sense intuitively. If we ask for
too much accuracy relative to the system size
($\alpha > 1$ in (\ref{eq:accsys}))
then tau-leaping's built-in bias outweighs its cheapness,
or, equivalently, the required stepsize is so small that
tau-leaping works harder than exact simulation.
The remainder of the table will be considered in section \ref{sec:CA}.

\begin{table}
\begin{center}
\begin{tabular}{|llll|}
  \hline
 Monte Carlo method    &  Computational complexity & unbiased? & Most efficient\\
  \hline \hline
  MC + exact simulation  &
          $\Theta(N^{2\alpha}+N)$ & Yes & Never\\[0.75ex]
 MC + tau-leaping &
          $\Theta(N^{3\alpha-1} +N^\alpha)$ & No & Never \\[0.75ex]
 MC + midpt.~or trap.~tau-leap &
          $\Theta(N^{2.5\alpha-1}  + N^{\alpha/2})$ & No & $\frac12 < \alpha \le \frac23$\\[0.75ex]
MC + Euler for diff.~approx.  & $\Theta(N^{3\alpha-1} + N^\alpha)$ & No & Never \\[0.75ex]
 MLMC $+$ E-M for  diff.~approx. &$\Theta(N^{2\alpha-1} + N^\alpha) $ & No & $\alpha \ge \frac23$ \\[0.75ex]
          biased MLMC  tau-leaping & $\Theta( N^{2\alpha -1}(\log N)^2  + N^\alpha)$ & No  & $\alpha\ge \frac23$ \\[0.75ex]
          unbiased MLMC tau-leaping & $\Theta(N^{2\alpha-1}(\log N)^2+N)$  & Yes & $\alpha \ge 1$ \\[0.75ex]
 \hline
 \end{tabular}
 \caption{Computational cost for different Monte Carlo methods, as $N\to \infty$.  The final column indicates when each method is most efficient, in terms of the parameter $\alpha$, up to factors involving logarithms.
 }
\label{tabNEW}
\label{tab1}
\end{center}
\end{table}

 We also mention that a crude and inexpensive approximation to the required expected value
can be computed by simply simulating the deterministic mass action ODE approximation to
(\ref{eq:model_N0}), which is often referred to as the reaction rate equation \cite{AndKurtz2011,AK2015}.  Depending upon the choice of functional $f$ and the underlying model \eqref{eq:model_family}, the bias from the ODE approximation can range from zero (in the case of a linear $\lambda_k$ and  linear function $f$), to order $N^{-1/2}$ (for example, when $f(X^{N}(\cdot)) = \sup_{t\le T} |X^N(t) - c(t)|$, where $c$ is the ODE approximation itself).  As we are interested in the fluctuations inherent to the stochastic model, we view $\alpha = \hhf$ as a natural cut-off in the relationship
    (\ref{eq:accsys}).

In addition to
the asymptotic complexity counts in
Table~\ref{tab1}, another important feature
of a method is
the availability of
computable
\emph{a posteriori}
confidence interval information.
As indicated in the table, two of the methods considered here,
exact simulation with Monte Carlo and an appropriately constructed multilevel tau-leaping, are
unbiased. The sample mean, accompanied
by an estimate of the overall variance, can then be delivered with a
computable
confidence interval.
By contrast, the remaining methods in the table are biased:
tau-leaping and Euler-Maruyama introduce discretization errors and the diffusion
approximation perturbs the underlying model. Although the asymptotic leading order
of these biases can be estimated, useful a posteriori upper bounds cannot be computed
straightforwardly in general,
making these approaches
much less attractive for
reliably achieving a target accuracy.

Based on the range of methods analysed here in an asymptotic regime
that couples system size and target accuracy,
three key messages are
\begin{itemize}
\item simulating exact samples alone is never advantageous,
 \item  even assuming there is no bias to the underlying model, simulating at the level of the the diffusion approximation is only marginally advantageous,
 \item tau-leaping can offer advantages over exact simulation, and an appropriately designed
     version of
     multilevel tau-leaping (which combines exact and tau-leaped samples)
     offers an unbiased method that is efficient over a wide range of
     accuracy requirements.
\end{itemize}

\section{Approximation methods}
\label{sec:approximations}

In this section, we briefly review the two alternatives to exact simulation of \eqref{eq:model_family} we study in this paper: tau-leaping and an Euler-Maruyama discretization of a diffusion approximation.

\subsection{Tau-Leaping}
\label{subsec:tau}
Tau-leaping
\cite{Gill2001}
is a computational method
that generates Euler-style approximate paths for the continuous-time  Markov chain
\eqref{eq:model_family}.
The basic idea is
to hold the intensity functions fixed over a time interval $[t_n,
t_n+h]$ at the values $\lambda_k(X^N(t_n))$, where $X^N(t_n)$ is the
 state of the system at time $t_n$, and, under this simplification, compute the
number of times each reaction takes place over this period.
 As the
waiting times for the reactions are exponentially distributed, this
leads to the following algorithm, which simulates up to a time of $T>0$.
For $x \ge 0$ we will write
Poisson$(x)$ to denote a sample from the Poisson
distribution with parameter $x$, with all
such samples
being independent of each other and of all
other sources of randomness used.

\begin{algorithm}[Euler tau-leaping]   \label{alg:Euler}
 Fix $h>0$.  Set $Z^N_h (0) = x_0$, $t_0 = 0$,
  $n=0$ and repeat the following until $t_{n} = T$:
  \begin{enumerate}[$(i)$]
   \item Set $t_{n+1} = t_{n} + h$.  If $t_{n+1} \ge T$, set $t_{n+1} = T$ and $h = T - t_n$.
  \item For each $k$, let $\Lambda_k = \text{Poisson}(\lambda_k(Z^N_h(t_n))h)$.
  \item Set $Z^N_h(t_{n+1}) = Z^N_h(t_{n}) + \sum_k \Lambda_k \zeta_k$.
   \item Set $n \leftarrow n+1$.
  \end{enumerate}
\end{algorithm}

Analogously
to \eqref{eq:model_family}, a path-wise representation of Euler tau-leaping defined for all $t\ge 0$
can be given through a random time change of Poisson processes:
\begin{equation}
  Z^N_h(t) = Z^N_h(0) + \sum_k \frac1N Y_k \left( N \int_0^t \lambda_k(Z^N_h (\eta_h(s)))
    ds  \right)\zeta_k,
  \label{eq:RTC_tau}
\end{equation}
where the $Y_k$ are as before, and $\displaystyle \eta_h(s) \eqdef \left \lfloor \frac{s}{h} \right \rfloor h$. Thus, $Z^N_h ( \eta_h(s) ) = Z^N_{h}(t_n)$ if $t_n\le s < t_{n+1}$. 
As the values of $Z_h^N$ can  go negative, the functions $\lambda_k$ must be defined outside of $\Z^d_{\ge 0}$.  One option is to simply define  $\lambda_k(x) = 0$ for $x\notin \Z^d_{\ge0}$, though other options exist \cite{Anderson2007b}.

\subsection{Diffusion approximation}
\label{subsec:diffusion}

The tau-leaping algorithm utilizes a time-stepping method to directly approximate the underlying model \eqref{eq:model_family}.  Alternatively, a diffusion approximation arises by perturbing the underlying model into one which can be discretized more efficiently.

Define the function $F$ via
\[
	F(x) = \sum_k  \lambda_k(x) \zeta_k.
\]
By the functional central limit theorem,
\begin{equation}\label{eq:another_BM_approx}
  \frac{1}{\sqrt{N}}\left[ Y_k(N u) - Nu\right] \approx W_k(u),
  \end{equation}
where $W_k$ is a standard Brownian motion.  Applying \eqref{eq:another_BM_approx} to \eqref{eq:model_family} yields
 \begin{equation*}
   	X^N(t) \approx X^N(0) +  \int_0^t F(X^N(s))ds + \sum_k \frac{1}{\sqrt{N}} W_k\left( \int_0^t  \lambda_k(X^N(s)) ds\right) \zeta_k,
   \end{equation*}
   where the $W_k$ are independent standard Brownian motions.  This implies that
 $X^N$ can be
 approximated by the process $D^N$ satisfying
   \begin{equation}\label{eq:4567898}
   	D^N(t) =  D^N(0) +  \int_0^t F(D^N(s))ds + \sum_k \frac{1}{\sqrt{N}} W_k\left( \int_0^t  \lambda_k(D^N(s)) ds\right) \zeta_k,
   \end{equation}
   where $D^N(0) = X^N(0)$.
An equivalent, and more prevalent, way to represent $D^N$ is via the It\^o representation
\begin{equation}
	D^N(t) =  D^N(0) +  \int_0^t F(D^N(s))ds + \sum_k \frac{1}{\sqrt{N}} \zeta_k  \int_0^t \sqrt{ \lambda_k(D^N(s)) }dW_k(s),
\label{eq:DNint}
\end{equation}
which is often written in the differential form
\begin{equation}\label{eq:ItoSDE}
	dD^N(t) = F(D^N(t))dt + \sum_k \frac{1}{\sqrt{N}} \zeta_k \sqrt{ \lambda_k(D^N(t))}dW_k(t),
\end{equation}
where  the $W_k$ of \eqref{eq:ItoSDE} are not necessarily the same as those in \eqref{eq:4567898}.

The SDE system
(\ref{eq:ItoSDE})
 is known as a {\em Langevin} approximation in the biology and chemistry literature, and a {\em diffusion} approximation in probability \cite{AK2015,Wilkinson2011}.
We note the following points.
\begin{itemize}
\item
The diffusion coefficient, often termed the ``noise'' in the system, is $\Theta(\frac{1}{\sqrt{N}})$, and hence,
in our setting is small relative to the drift.
\item
The  diffusion coefficient involves square roots. Hence, it is critical that the intensity functions $\lambda_k$ only take values in $\R_{\ge 0}$ on the domain of the solution.  This is of particular importance in the population process setting where the solutions of the underlying model \eqref{eq:model_family} naturally satisfy a non-negativity constraint whereas the SDE solution paths cannot be guaranteed to
remain non-negative in general.  In this case
one reasonable representation, of many,  would be
\begin{equation}\label{eq:diffusion_pos}
dD^N(t) = F(D^N(t))dt + \sum_k \frac{1}{\sqrt{N}} \zeta_k \sqrt{[ \lambda_k(D^N(s))]^+ }dW_k(s),
\end{equation}
where $[x]^+ = \max\{x,0\}$.  Another reasonable option would be to use a process with reflection \cite{LW2016}.
\item
The coefficients of the SDE are not globally Lipschitz in general,
and hence standard convergence theory
for numerical methods, such as that in
\cite{KloedenPlaten92},
is not applicable.
Examples of nonlinear SDEs for which standard Monte Carlo and
multilevel Monte Carlo,
when combined with and
Euler-Maruyama discretization with a uniform timestep, fail to
produce a convergent algorithm
have been pointed out in the literature
\cite{HJK11,HJK11b}.
The question of which classes of reaction systems
lead to well-defined SDEs and
which
discretizations converge at the traditional rate therefore remains open.
\end{itemize}

  In this work, to get a feel for the best possible computational complexity
that can arise from the Langevin approximation, we will study the case where the bias that arises from switching models from $X^N$ to $D^N$ is zero. We will also assume that, even though the diffusion coefficients involve square roots and are therefore not generally globally Lipschitz, the Euler-Maruyama method has a bias of order $\Theta(h)$.  
We will find that even in this idealized light, the asymptotic computational complexity of Euler-Maruyama on a diffusion approximation combined with either a standard or a multilevel implementation is only marginally better than the corresponding  computational complexity bounds for multilevel tau-leaping. In particular, they differ only in a  factor of a logarithm of the scaling parameter.

Finally, due to the fact that the diffusion approximation itself already has a difficult to quantify bias, we will not consider higher order methods \cite{AndMattTrap}, or even unbiased methods \cite{GR2005}, for this process.

\section{Complexity analysis}
\label{sec:CA}

In this section we establish the results given in Table~\ref{tab1}.
In subsection \ref{sec:MC},  we derive  the first four rows, whereas in subsection \ref{sec:MLMC} we discuss the multilevel framework and establish rows five, six, and seven.

\subsection{Complexity analysis of standard Monte Carlo approaches}
\label{sec:MC}

\subsubsection{Exact Sampling and Monte Carlo}
\label{subsec:exactMC}

 By Assumption \ref{assump:costperpath} the expected number of system updates
required to generate a single exact sample path is  $\Theta(N)$.  Letting
\[
	\delta_{N} = \textsf{Var}(f(X^N(\cdot))),
\]
in order to get a standard deviation below $\varepsilon_N$ we require 
\[
	n^{-1} \delta_{N} \le  \varepsilon_N^2 \implies n \ge \delta_N\varepsilon_N^{-2}+1.
\]
Thus, the total computational complexity of making the desired approximation is
\[
	\Theta(nN)= \Theta(\delta_N \varepsilon_N^{-2} N+N) = \Theta(\delta_N N^{2\alpha+1} +N) .
\]
By Assumption \ref{assump:variance}, $\delta_N = \Theta(N^{-1})$, yielding an overall complexity of
$
 \Theta(N^{2\alpha}+N)$, as given in the first row of Table \ref{tab1}.

\subsubsection{Tau-leaping and Monte Carlo}
\label{subsec:tauMC}

Suppose now that we use $ n $ paths of the tau-leaping
process
(\ref{eq:RTC_tau})
to construct the Monte Carlo estimator
$\hat \mu_n$
for $\E [ f(X^N(\cdot))]$.  By assumption \ref{assump:bias}, the bias is $\Theta(h)$, so we constrain ourselves to $h = \varepsilon_N$.  
Letting
\[
	\delta_{N,h} = \textsf{Var}(f(Z^N_h(\cdot)))
\]
we again require $n \ge \delta_{N,h}\varepsilon_N^{-2}+1$ to control the statistical error.  Since by Assumption \ref{assump:costperpath} there are $\Theta(h^{-1})$ expected operations per path generation, the total computational complexity for making the desired approximation is 
\[
	\Theta(nh^{-1}) = \Theta(\delta_{N,h}\varepsilon_N^{-3}+\varepsilon_N^{-1}) .
\]
By Assumption \ref{assump:variance}, $\textsf{Var} (f(Z^N_{h,i}(\cdot))) = \Theta(N^{-1})$,
giving an overall complexity of
$
\Theta(N^{3\alpha - 1}+N^{\alpha})
$, as reported in the second row of Table \ref{tab1}.

Weakly second order extensions to the tau-leaping
method
can lower the computational
complexity dramatically.
For example,
if we use the midpoint tau-leaping
process ${\cal Z}_h^N$ from
  \cite{AndersonGangulyKurtz}, by Assumption \ref{assump:bias} we can set $h = \sqrt{\varepsilon_N}$ and still achieve a bias of $\Theta(\varepsilon_N)$.   
Since by Assumption \ref{assump:variance} we need
$n \ge N^{-1}\varepsilon_N^{-2}+1$ paths to control the standard deviation, the complexity is
\[
	\Theta(n\cdot h^{-1})   =\Theta(N^{-1}\varepsilon_N^{-2.5}+\varepsilon_N^{-1/2}) = \Theta(N^{2.5\alpha-1}+N^{\alpha/2}),
\]
   as stated in the third row of Table~\ref{tab1}.
The same conclusion can also be drawn for
the trapezoidal method in
\cite{AndMaso2011}.

If methods are developed that are higher order in a weak sense, then further improvements can be gained.  In general, if a method is developed that is weakly of order $\rho$, then we may set $h = \varepsilon_N^{1/\rho}$ to achieve a bias of $\Theta(\varepsilon_N)$.  Still supposing a per-path variance of $\Theta(N^{-1})$, we again choose $n \ge N^{-1}\varepsilon_N^{-2} + 1$ paths, and find a complexity of
\[
	\Theta(n\cdot h^{-1})   =\Theta(N^{-1}\varepsilon_N^{-(2+\frac1\rho)}+\varepsilon_N^{-1/\rho}) = \Theta(N^{(2+\frac1\rho)\alpha -1}+N^{\alpha/\rho}).
\]
For example, if a third order method is developed, i.e., $\rho =3$, then this method becomes optimal for $\frac12 \le \alpha \le \frac34$.  To the best of the authors' knowledge, no such methods have yet been designed.

\subsubsection{Diffusion approximation and Monte Carlo}
\label{subsec:diffMC}

Given Assumptions \ref{assump:costperpath}, \ref{assump:bias}, and \ref{assump:variance}, the complexity analysis for the diffusion approximation with Euler-Maruyama is exactly the same as for Euler tau-leaping.  Hence,  we can again give an overall complexity of
$
\Theta(N^{3\alpha-1}+N^{\alpha})
$, as reported in the fourth row of Table \ref{tab1}.

\subsection{Multilevel Monte Carlo and complexity analysis}
\label{sec:MLMC}
In this section we study multilevel Monte Carlo
approaches and derive the results summarized in rows
five, six, and seven of Table~\ref{tab1}.

\subsubsection{Multilevel Monte Carlo and Diffusion Approximation}
\label{subsec:diffMLMC}

Here we specify and analyze an Euler-based multilevel method for the diffusion approximation, following the original framework of Giles \cite{Giles2008}.  

For some fixed $M>1$ we let $h_{\ell} = T \cdot M^{-\ell}$ for $\ell \in \{0,\dots,L\}$, where $T>0$ is a fixed terminal time.  Reasonable choices for $M$ include $M\in \{2,3,4,\dots , 7\}$, and $L$ is determined below.
Let $\DN_{h_\ell}$ denote the approximate process generated by Euler-Maruyama
applied to (\ref{eq:ItoSDE})
with a step size of $h_{\ell}$.
By Assumption \ref{assump:bias} 
 we may set $h_L = \varepsilon_N$, giving $L = \Theta( |\log \varepsilon_N| )$, so that the finest level
achieves the required order of magnitude for the bias.

Noting that
\begin{equation}\label{eq:490857}
	\E[f(D^N_{h_L}(\cdot))] = \E[f(D^N_{h_0}(\cdot))] +  \sum_{\ell = 1}^L \E[ f(D^N_{h_{\ell}}(\cdot)) - f(D^N_{h_{\ell-1}}(\cdot))],
\end{equation}
we use $i$ as an index over sample paths and let
\begin{align*}
	\QNl_{0} &\eqdef \frac{1}{n_{0}} \sum_{i = 1}^{n_{0}} f(\DN_{h_0,[i]}(\cdot)), \quad \text{and} \quad
	\QNl_{\ell} \eqdef \frac{1}{n_{\ell}} \sum_{i = 1}^{n_{\ell}} ( f(\DN_{h_\ell,[i]}(\cdot)) - f(\DN_{h_\ell-1,[i]}(\cdot))),
\end{align*}
for $\ell = 1, \dots, L$, where $n_0$ and the different $n_\ell$ have yet to be determined.  Note that the form of the estimator $\QNl_\ell$ above implies that the processes $\DN_{h_\ell}$ and $\DN_{h_{\ell-1}}$ will be \textit{coupled}, or constructed on the same probability space.  We consider here the case when  $(D^N_{h_\ell},D^N_{h_{\ell-1}})$ are coupled in the usual way by using the same Brownian path in the generation of each of the marginal processes.
Our (biased) estimator is then
\[
  \QNl \ \eqdef \
  \QNl_0 +
  \sum_{\ell = 1}^L \QNl_{\ell}.
\]
 Set
\[
	\delta_{N,\ell} = \textsf{Var}  ( f(\DN_{h_\ell}(\cdot)) - f(\DN_{h_\ell-1}(\cdot))).
\]
By Assumption \ref{assump:variance}, 
$\delta_{N,\ell} = \Theta(N^{-1}h_\ell^2+N^{-2}h_\ell)$ and 
$\textsf{Var}(f(D_{h_0}^N(\cdot))) = \Theta(N^{-1}).$ 
In  \cite{AHS2016} it is shown that under these circumstances, the computational complexity required is $\Theta(\varepsilon_N^{-2} N^{-1} + \varepsilon_N^{-1})$.
In the regime 
(\ref{eq:accsys})
this translates to 
$\Theta(N^{2\alpha - 1} + N^\alpha),$ as reported in the fifth row of Table \ref{tab1}.

\subsubsection{Multilevel Monte Carlo and tau-leaping}
\label{subsec:tauMLMC}

The use of multilevel Monte Carlo with tau-leaping for
continuous-time Markov chains of the form considered here
was proposed in
\cite{AndHigh2012},
where effective algorithms were devised.   Complexity results were
given in a non-asymptotic multi-scale
setting, with followup results in \cite{AHS2014}.
Our aim here is to customize
the approach in the
scaling regime
 (\ref{eq:accsys})
and thereby develop
easily
interpretable complexity bounds
that allow straightforward comparison with
other methods.
In this section $Z^N_{h_\ell}$ denotes a tau-leaping process generated with a step-size of $h_\ell = T\cdot M^{-\ell}$, for $\ell \in \{0,\dots,L\}$.

A major step in
\cite{AndHigh2012}
was to show that
a coupling technique  used for analytical purposes
in
\cite{AndersonGangulyKurtz,Kurtz82}
can also form the basis of a practical simulation algorithm.
Letting
$Y_{k,i},\ i \in \{1,2,3\}$, denote independent,
unit rate Poisson processes, we couple the exact and approximate tau-leaping processes in the following way,
\begin{align}
\begin{split}
	X^N(t) = &X^N(0) + \sum_{k}\frac{1}{N} Y_{k,1}\left( N \int_0^t  \lambda_k(X^N(s))  \wedge  \lambda_k(Z^N_{h_L} ( \eta_{L}(s)))  ds  \right)\zeta_k\\
	&\hspace{.1in} + \sum_{k} \frac{1}{N}Y_{k,2}\left( N \int_0^t [\lambda_k(X^N(s))  -  \lambda_k(X^N(s))  \wedge   \lambda_k(Z^N_{h_L} (\eta_{L}(s)))  ] ds \right)\zeta_k,
    \label{eq:Z_X1}
\end{split}\\
\begin{split}
	Z^N_{h_L} (t) = &Z^N_{h_L}(0) + \sum_{k} \frac{1}{N}Y_{k,1}\left(N  \int_0^t \lambda_k(X^N(s))  \wedge   \lambda_k(Z^N_{h_L} (\eta_{L}(s)))  ds \right)\zeta_k \\
	&\hspace{.1in} + \sum_{k} \frac{1}{N} Y_{k,3}\left(N  \int_0^t [ \lambda_k(
Z^N_{h_L} ( \eta_{L}(s)))  - \lambda_k(X^N(s))  \wedge  \lambda_k(
Z^N_{h_L} (\eta_{L}(s))) ] ds  \right)\zeta_k,
      \label{eq:Z_X2}
    \end{split}	
\end{align}
where $a\wedge b$ denotes $\min\{a,b\}$ and $\eta_L(s) = \lfloor s/h_L\rfloor h_L$.
Sample paths of
    (\ref{eq:Z_X1})--(\ref{eq:Z_X2})
can be generated with a natural extension of the
next reaction method
 or Gillespie's algorithm, see \cite{AndHigh2012}, and for $h_L \ge N^{-1}$ the
complexity required for the generation of a realization $(X^N,Z_{h_{L}}^N)$ remains at the
$\Theta(N)$ level.
	The coupling of two approximate processes, $Z^N_{h_{\ell}}$ and $Z^N_{h_{\ell-1}}$, takes the similar form
\begin{align}
\begin{split}
	Z^N_{h_{\ell}}&(t)  = Z^N_{h_{\ell}}(0) + \sum_{k} \frac{1}{N} Y_{k,1}\left(N  \int_0^t \lambda_k(Z^N_{h_{\ell}} ( \eta_{\ell}(s))) \wedge \lambda_k(Z^N_{h_{\ell-1}} ( \eta_{\ell-1}(s))) ds \right)\zeta_k \\
	& +  \sum_k \frac{1}{N} Y_{k,2}\left(N \int_0^t [ \lambda_k(Z^N_{h_{\ell}} (\eta_{\ell}(s))) -  \lambda_k(Z^N_{h_{\ell}} (\eta_{\ell}(s))) \wedge \lambda_k(Z^N_{h_{\ell-1}} (\eta_{\ell-1}(s))) ] ds \right)\zeta_k,  \label{eq:Z1}
    \end{split}\\
    \begin{split}
	Z^N_{h_{\ell-1}}&(t) = Z^N_{h_{\ell-1}}(0) + \sum_{k} \frac{1}{N}  Y_{k,1}\left(N   \int_0^t \lambda_k(Z^N_{h_{\ell}} ( \eta_{\ell}(s))) \wedge \lambda_k(Z^N_{h_{\ell-1}} ( \eta_{\ell-1}(s)))  ds  \right)\zeta_k\\
	& +  \sum_k \frac{1}{N}Y_{k,3}\left( N \int_0^t [ \lambda_k(Z^N_{h_{\ell-1}} ( \eta_{\ell-1}(s))) -  \lambda_k(Z^N_{h_{\ell}} ( \eta_{\ell}(s))) \wedge \lambda_k(Z^N_{h_{\ell-1}} (\eta_{\ell-1}(s))) ] ds \right)\zeta_k,
    \label{eq:Z2}
    \end{split}	
\end{align}
where
$\eta_{\ell}(s) \eqdef \lfloor s/h_{\ell}\rfloor h_{\ell}$.
The  pair
(\ref{eq:Z1})--(\ref{eq:Z2}) can be sampled
at the same $\Theta(h_{\ell}^{-1})$ cost as a single
tau-leaping path, see
\cite{AndHigh2012}.

For  $L$ as yet to be determined, and  noting the identity
\begin{equation}
    \EE[ f( X^N(\cdot)) ]
         =
    \EE[ f( X^N(\cdot)) -  f ( Z^N_L(\cdot))]
     +
   \sum_{\ell = {1}}^{L}
    \EE[ f( Z^N_{h_{\ell}}(\cdot) ) - f ( Z^N_{h_{\ell-1}}(\cdot) ) ]
     +
    \EE[ f(Z^N_{h_{0}}(\cdot)) ],
 \label{eq:tel2}
\end{equation}
we define estimators
for the three terms above via
\begin{align}
\begin{split}
	\QNl_{E} &\eqdef \frac{1}{n_{E}} \sum_{i = 1}^{n_{E}} (f(X^N_{[i]}(\cdot)) - f(Z^N_{h_L,[i]}(\cdot))) ,\\
	\QNl_{\ell} &\eqdef \frac{1}{n_{\ell}} \sum_{i = 1}^{n_{\ell}} ( f(Z^N_{h_{\ell},[i]}(\cdot)) - f(Z^N_{h_{\ell} - 1,[i]}(\cdot))), \quad \text{for } \ell \in \{1,\dots,L\}, \\
	\QNl_{0} &\eqdef \frac{1}{n_{0}} \sum_{i = 1}^{n_{0}} f(Z^N_{h_{0},[i]}(\cdot)),
	\end{split}
	\label{eq:unbiased_estimators}
\end{align}
so that
\begin{equation}
	\QNl  \eqdef  \QNl_E +  \sum_{\ell = 1}^L \QNl_{\ell} + \QNl_0
	\label{eq:unbiased_MLMC}
\end{equation}
is an unbiased estimator for $\E [f(X^N(\cdot))]$.
Here,
$\QNl_{E}$ uses the coupling
    (\ref{eq:Z_X1})--(\ref{eq:Z_X2})
between exact paths and
tau-leaped paths of stepsize $h_L$,
$\QNl_{\ell}$ uses the coupling
    (\ref{eq:Z1})--(\ref{eq:Z2})
between tau-leaped paths of stepsizes $h_{\ell}$  and
$h_{\ell-1}$, and
$\QNl_{0}$ involves single tau-leaped paths of
stepize $h_{0}$.
Note that the algorithm implicit in \eqref{eq:unbiased_MLMC} produces an unbiased estimator, whereas the estimator is biased if  $\QNl_E$ is left off, as will sometimes be desirable.  Hence, we will refer to estimator $\QNl$ in \eqref{eq:unbiased_MLMC} as the unbiased estimator, and will refer to
\begin{equation}\label{eq:biased_version}
	\QNl_B \eqdef \sum_{\ell = 1}^L \QNl_{\ell} + \QNl_0
\end{equation}
as the biased estimator.
For both the biased and unbiased estimators, the number of paths at each level,
$n_0$, $n_{\ell}$ and $n_{E}$,
will be chosen to ensure an overall estimator standard deviation 
of $\varepsilon_N$.

We consider the biased and unbiased versions of tau-leaping multilevel Monte Carlo separately.

\vspace{.125in}

\noindent \textbf{Biased multilevel Monte Carlo tau-leaping}

\vspace{.125in}

Here we consider the estimator $\QNl_B$ defined in  \eqref{eq:biased_version}. By Assumption \ref{assump:bias}  $| \E[f(X^N(\cdot))] - \E[f(Z^N_{h_L}(\cdot))] | = \Theta(h_L)$.  Hence, in order to control the bias we begin by choosing $h_L = \varepsilon_N$ and so $L = \Theta(\log(1/\varepsilon_N)) = \Theta(\log N)$.

For $\ell \in\{1,\dots,L\}$, let $C_\ell$  be the expected number
of random variables required to generate a single pair of coupled trajectories at level $\ell$ and let $\delta_{N,\ell}$ be the variance of the relevant processes on level $\ell$.  Let $C_0$  be the expected number of random variables  required to generate a single trajectory at the coarsest level. To find $n_\ell$, $\ell \in \{0,\dots,L\}$,  we solve the following optimization problem, which ensures that the variance of $\QNl_B$ is no greater than  $\varepsilon_N^2$:
 \begin{align}
&\underset{n_{\ell}}{\text{minimize}} \hspace{.15in} \sum_{\ell = 0}^L n_{\ell}C_{\ell},\label{eq:LM1}\\
&\text{subject to} \hspace{.1in} \sum_{\ell = 0}^L \frac{\delta_{N,\ell} }{n_{\ell}}= \varepsilon_N^2.\label{eq:LM2}
\end{align}
     We use  Lagrange multipliers.  Since we have $C_{\ell}=K\cdot h_{\ell}^{-1}$, for some fixed constant $K$, the optimization problem above is solved at solutions to
\begin{equation*}
 \nabla_{n_0,\dots,n_L,\lambda}\left(
 \sum_{\ell = 0}^L n_{\ell}K\cdot h_{\ell}^{-1}+\lambda\left(\sum_{\ell = 0}^L \frac{\delta_{N,\ell} }{n_{\ell}}- \varepsilon_N^2\right)\right) = 0.
\end{equation*}
Taking derivatives with respect to $n_{\ell}$  and setting each derivative to zero yields,
\begin{equation}
\label{eq:nl}
n_{\ell}= \sqrt{\tfrac{\lambda}{K}\delta_{N,\ell}h_{\ell}}, \qquad \text{ for } \ell \in\{0,1,2,\dots,L\}
\end{equation}
for some $\lambda\ge 0$. Plugging \eqref{eq:nl} into \eqref{eq:LM2} gives us,
\begin{equation}\label{eq:78966798}
\sum_{\ell = 0}^L\sqrt{ \frac{\delta_{N,\ell}}{h_{\ell}}}= \sqrt{\tfrac{\lambda}{K}}\cdot \varepsilon_N^2
\end{equation}
and hence by Assumption \ref{assump:variance}
\begin{equation}\label{eq:09877890}
 \sqrt{\tfrac{\lambda}{K}}=\sum_{\ell = 0}^L\sqrt{ \tfrac{\delta_{N,\ell}}{h_{\ell}}} \le CL\varepsilon_N^{-2}N^{-1/2},
\end{equation}
where $C$ is a constant.
Noting that $L=\Theta(\log(\varepsilon_N^{-1}))$, we have
\[
\tfrac{\lambda}{K}=\Theta\left(\varepsilon_N^{-4}\left(\log \varepsilon_N\right)^2N^{-1}\right).
\]
Plugging this back into \eqref{eq:nl}, and recognizing that at least one path must be generated to achieve the desired accuracy, we find
\[
	n_{\ell} = \Theta(\varepsilon_N^{-2} N^{-1}h_\ell L+ 1).
\]
Hence, the overall computational complexity is
\begin{align*}
\sum_{\ell=0}^L n_\ell K h_\ell^{-1} &= \Theta \left( \sum_{\ell=0}^L\varepsilon_N^{-2} N^{-1}h_\ell L h_\ell^{-1}  + \sum_{\ell=0}^L h_{\ell}^{-1}\right) = \Theta \left( \varepsilon_N^{-2} N^{-1}(\log \varepsilon_N)^2+ \varepsilon_N^{-1}\right)\\
&= \Theta \left( N^{2\alpha-1}(\log N)^2+ N^{\alpha}\right),
\end{align*}
recovering row six of Table~\ref{tab1}.

Note that the  computational complexity reported for this biased version of multilevel Monte Carlo tau-leaping is, up to logarithms,  the same as that for multilevel Monte Carlo on the diffusion approximation.  However, none of the generous assumptions we made for the diffusion approximation were required.

\vspace{.125in}

\noindent \textbf{Unbiased multilevel Monte Carlo tau-leaping }
\vspace{.125in}

The first observation to make is that the telescoping sum \eqref{eq:tel2} implies that the method which utilizes  $\E[f(X^N(\cdot)) - f(Z^N_{h_L}(\cdot))]$ at the finest level is unbiased for \textit{any} choice of $h_L$.  That is, we are no longer constrained to choose $L = \Theta(|\log \varepsilon_N|)$.

Assume that $h_L \ge N^{-1}$. Let $C_E$ be the expected number of random variables required to generate a single pair of the coupled exact and tau-leaped processes when the tau-leap discretization is $h_L$. 
To determine $n_\ell$ and $n_E$, we still solve an  optimization problem,

 \begin{align}
&\underset{n_{\ell}}{\text{minimize}} \hspace{.15in} \sum_{\ell = 0}^L n_{\ell}C_{\ell} + n_LC_E,\label{eq:UM1}\\
&\text{subject to} \hspace{.1in} \sum_{\ell = 0}^L \frac{\delta_{N,\ell} }{n_{\ell}}  + \frac{\delta_{N,E}}{n_E}= \varepsilon_N^2,\label{eq:UM2}
\end{align}
where $C_\ell$ and $\delta_{N,\ell}$ are as before and $\delta_{N,E} = \textsf{Var}(f(X^N(\cdot)) - f(Z^N_{h_L}(\cdot)))$.

Using Lagrange multipliers again, we obtain,
\begin{equation}
\label{eq:nll}
	n_{\ell} = \sqrt{\frac{\lambda\delta_{N,  \ell}}{C_{\ell}}}   \qquad \text{ for } \ell \in\{0,1,2,\dots,L\}
\end{equation}
and 
\begin{equation}
\label{eq:ne}
	n_{E} = \sqrt{\frac{\lambda\delta_{N,  E}}{C_E}}.
\end{equation}
Plugging back into \eqref{eq:UM2} and noting that by, Assumption  \ref{assump:costperpath}, $C_{\ell} = \Theta(h_{\ell}^{-1})$ and  $C_{E} = \Theta(N)$ yields
\begin{equation}\label{eq:5676567}
 \sqrt{\lambda}=\varepsilon_N^{-2}\left(\sum_{\ell = 0}^L\sqrt{ \delta_{N,\ell}C_{\ell}} +\sqrt{ \delta_{N,E}C_{E}}\right)\le C(L\varepsilon_N^{-2}N^{-1/2} +\varepsilon_N^{-2}\sqrt{h_{L}} ).
\end{equation}
Therefore,  plugging \eqref{eq:5676567} back into \eqref{eq:nll} an \eqref{eq:ne} and noting $n_{\ell}  \ge 1$ and $n_E \ge 1$, we get
\begin{equation*}
	n_{\ell} =\sqrt{\frac{\lambda\delta_{N,  \ell}}{C_{\ell}}}  + 1=O\left(\left(L\varepsilon_N^{-2}N^{-1} +\varepsilon_N^{-2}\sqrt{\frac{h_{L}}{N}}\right) h_{\ell}+ 1\right)    \qquad \text{ for } \ell \in\{0,1,2,\dots,L\}
\end{equation*}
and 
\begin{equation}
	n_{E} = \sqrt{\frac{\lambda\delta_{N,  \ell}}{C_{\ell}}} + 1 =  O(L\varepsilon_N^{-2}N^{-3/2}h_{L} ^{1/2}+\varepsilon_N^{-2}N^{-1}h_{L}+ 1).
\end{equation}
As a result the total complexity is
\begin{align*}
g(h_L) &= O(\varepsilon_N^{-2}N^{-1}L ^2+\varepsilon_N^{-2}\sqrt{\frac{h_{L}}{N}}L+ h_L^{-1}+ \varepsilon_N^{-2}\sqrt{\frac{h_{L}}{N}}L+\varepsilon_N^{-2}h_{L} + N)\\
&\le O(\varepsilon_N^{-2}N^{-1}L ^2+2\varepsilon_N^{-2}\sqrt{\frac{h_{L}}{N}}L+\varepsilon_N^{-2}h_{L} + 2N)\tag{since $h_L^{-1} \le N$}\\
&=O(2\varepsilon_N^{-2}N^{-1}L ^2+2\varepsilon_N^{-2}h_{L} + 2N)\tag{using that $2ab \le a^2 + b^2$}.
\end{align*}
 It is relatively easy to show that the last line above is minimized at
\begin{equation}\label{eq:89769876}
h_L = \frac{2}{(\log 2)^2 N} \text{LambertW}\left(\frac{N}{2/(\log 2)^2}\right) \approx \frac{2}{(\log 2)^2 N} \log\left( \frac{N}{2/(\log 2)^2)}\right).
\end{equation}

Hence, taking $h_L = \Theta(N^{-1} \log N)$, we  have $(\log h_L)^2 = \Theta((\log N)^2)$ and this method  achieves a total computational complexity of leading order
\[
	\Theta(\varepsilon_N^{-2}N^{-1}(\log N)^2+\varepsilon_N^{-2}N^{-1}\log N + N) = \Theta(\varepsilon_N^{-2}N^{-1}(\log N)^2+ N) = \Theta(N^{2\alpha-1}(\log N)^2+ N) ,
\]
as reported in the last  row of Table \ref{tab1}.

Note here that if we choose $h_L = \frac{1}{N}$ we get the same order of magnitude for the computational complexity.  However the $h_L$ in \eqref{eq:89769876} is the optimized solution, meaning the leading order constant should be better and we will see this in Figure \ref{fig:example_compare1} and Figure \ref{fig:example_compare54}  in the next section.

\section{Computational results}
\label{sec:comp}

In this section we provide numerical evidence for the sharpness of the computational complexity analyses provided in Table \ref{tab1}. We will measure complexity by total number of random variables utilized. We emphasize that these experiments use extreme parameter choices solely for the purpose
of testing the sharpness of the delicate asymptotic bounds. 

\begin{example}

We consider the classically scaled stochastic model for the following reaction network (see \cite{AK2015})
\[
	S_1+S_2\overset{k_1/N}{\underset{k_2}{\rightleftarrows}} S_3\overset{k_3}{\rightarrow}S_2+S_4.
\]
Letting $X_i(t)$ give the number of molecules of species $S_i$ at time $t$, and letting $X^N(t) = X(t)/N$, the stochastic equations are
\begin{align*}
	X^N(t) = X^N(0) + \frac{1}{N} &Y_1\left( N k_1\int_0^t X_1^N(s)X_2^N(s) ds\right)\left[\begin{array}{r}
	-1\\
	-1\\
	1\\
	0
	\end{array}\right] \\
	 + \frac{1}{N} &Y_2\left( N k_2 \int_0^t X^N_3(s) ds \right) \left[\begin{array}{r}
          1\\
	 1\\
	-1\\
	0
	\end{array}\right]\\
         + \frac{1}{N} &Y_3\left( N k_3 \int_0^t X^N_3(s) ds \right) \left[\begin{array}{r}
          0\\
	 1\\
	-1\\
	1
	\end{array}\right],
\end{align*}
where we assume $X^N(0)\to (0.2,0.2,0.2,0.2)^T$, as $N\to \infty$.  Note that the intensity function $\lambda_1(x)= \kappa_1 x_1 x_2$ is globally Lipschitz on the domain of interest as that domain is bounded (mass is conserved in this model).

We implemented different Monte Carlo simulation methods for the estimation of $\E [X_1^N(T)]$ to an accuracy of $\varepsilon_N= N^{-\alpha}$ for both $\alpha=1$ and $\alpha=5/4$.  Specifically, for each of the order one methods we chose a step size of  $h = \varepsilon_N$ and  required the variance of the estimator to be $\varepsilon_N^2$.  For midpoint tau-leaping, which has a weak order of two, we chose $h = \sqrt{\varepsilon_N}$.  For the unbiased multilevel Monte Carlo method we chose the finest time-step according to \eqref{eq:89769876}.    We do not provide estimates for Monte Carlo combined with exact simulation as those computations were too intensive to complete to the target accuracy.

For our numerical example we chose $T = 1$ and $X(0) = \lceil N\cdot  [0.2, 0.2, 0.2, 0.2]^T\rceil$ with $X^N(0)=X(0)/N$.  Finally, we chose  $k_1 = k_2 =  k_3 = 1$ as our rate constants. In Figure \ref{fig:example_simple}, we provide log-log plots of the computational complexity required to solve this problem for the different Monte Carlo methods to an accuracy of   $\varepsilon_N = N^{-1}$, for each of  
\[
	N\in\{2^{13},2^{14},2^{15},2^{16},2^{17}\}.
\] 
In Figure \ref{fig:example_54}, we provide log-log plots for the computational complexity required to solve this problem for the different methods to an accuracy of  $\varepsilon_N = N^{-\frac{5}{4}}$, for each of
\[
	N\in\{2^{9},2^{10},2^{11},2^{12}, 2^{13}\}.
\]

 \begin{figure}
\centering
     \includegraphics[width = \textwidth]{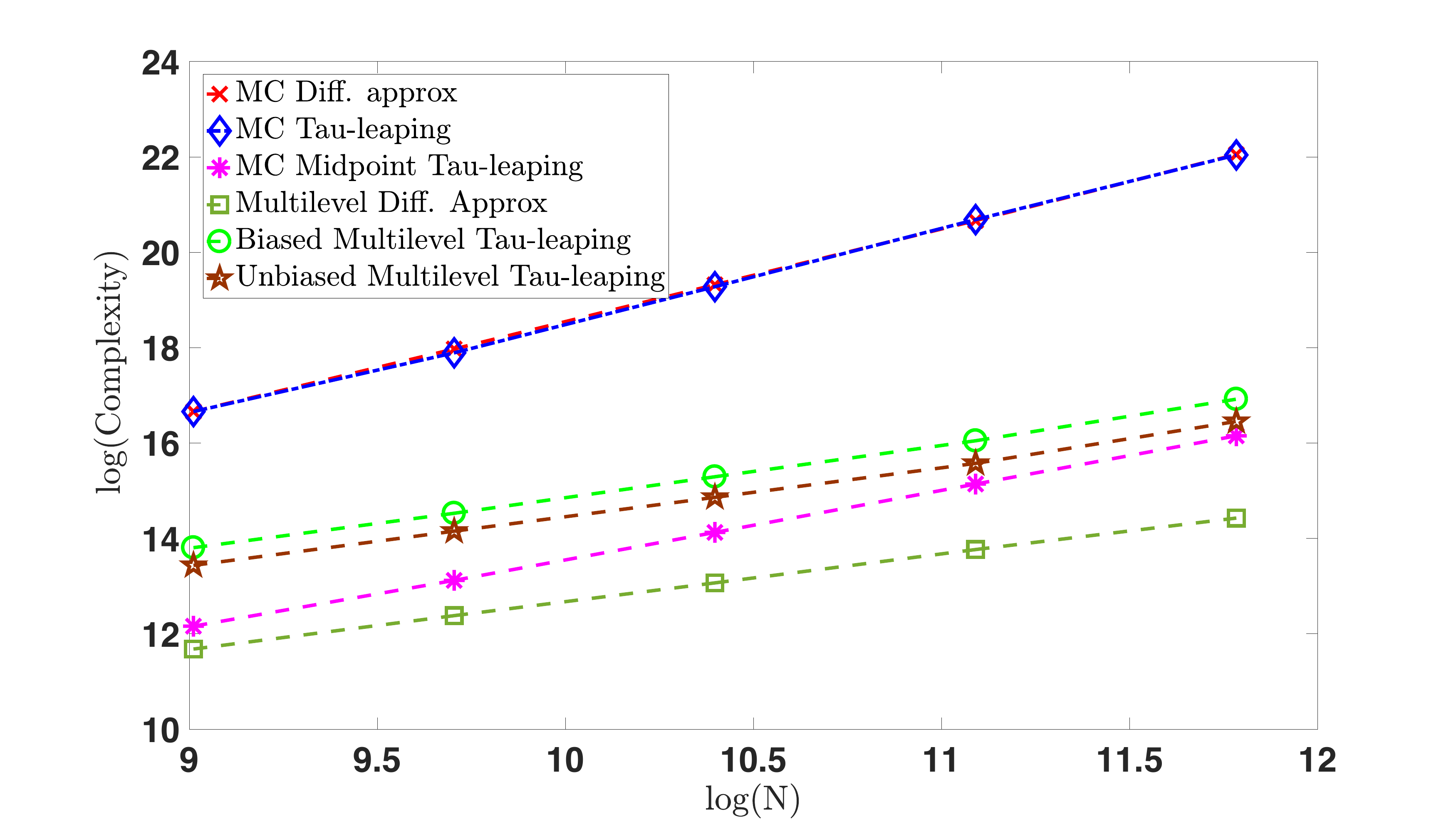}
	\caption{Log-log plots of the computational complexity for the different Monte Carlo methods with varying $N\in\{2^{13},2^{14},2^{15},2^{16},2^{17}\}$ and  $\varepsilon_N = N^{-1}$.}
	\label{fig:example_simple}
\end{figure}

\begin{figure}
\centering
     \includegraphics[width = \textwidth]{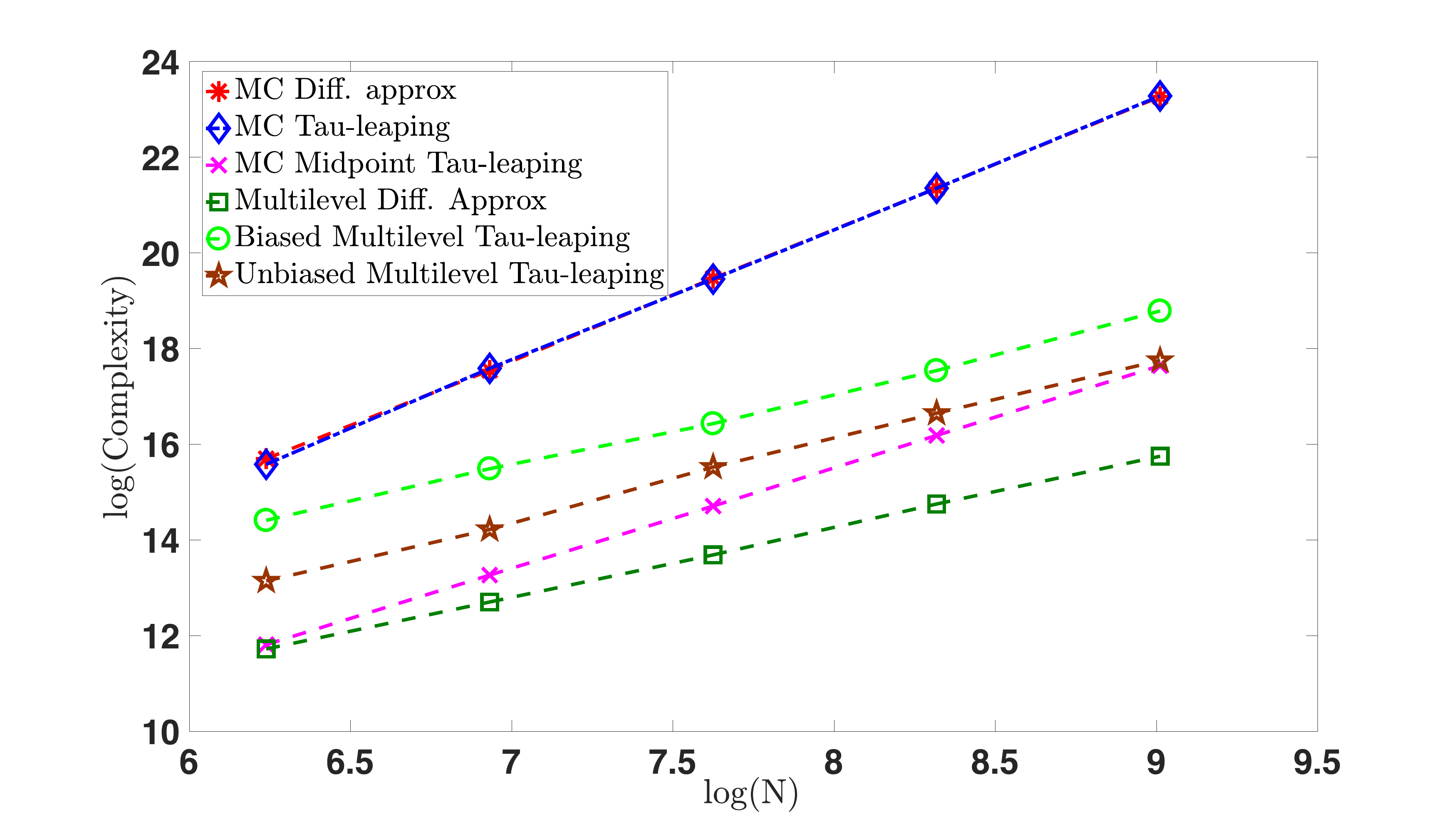}
	\caption{Log-log plots of the computational complexity for the different Monte Carlo methods with varying $N\in\{2^{9},2^{10},2^{11},2^{12}, 2^{13}\}$ and  $\varepsilon_N = N^{-\frac{5}{4}}$.}
\label{fig:example_54}
\end{figure}

Tables \ref{table1} and \ref{table2} provide  the  estimator standard deviations for the different Monte Carlo methods with $\varepsilon_N = N^{-1}$ and  $\varepsilon_N = N^{-\frac{5}{4}}$, respectively.  The top line provides the target standard deviations.

\begin{table}
\begin{center}
\begin{tabular}{|c|c|c|}\hline
Method  & estimator standard deviations \\ \hline
  & $2^{-13},2^{-14},2^{-15},2^{-16},2^{-17}$\\ \hline
MC and Diff. approx  & $2^{-13.10},2^{-14.02},2^{-15.02},2^{-16.01},2^{-17.00} $\\ \hline
MC and Tau-leaping  & $ 2^{-13.09},2^{-14.01},2^{-15.01},2^{-16.01},2^{-17.00} $\\ \hline
MC and Midpoint Tau-leaping  & $2^{-13.09},2^{-14.04},2^{-15.03},2^{-16.00},2^{-17.01}$ \\ \hline
Multilevel Diff. approx  & $2^{-13.20},2^{-14.15},2^{-15.11},2^{-16.09},2^{-17.07} $\\ \hline
Biased Multilevel Tau-leaping & $2^{-13.44},2^{-14.39},2^{-15.39},2^{-16.38},2^{-17.32} $\\ \hline
Unbiased Multilevel Tau-leaping & $2^{-13.29},2^{-14.28},2^{-15.26},2^{-16.21},2^{-17.18} $\\ \hline
\end{tabular}
\end{center}
\caption{Actual estimator standard deviations when $\varepsilon_N = N^{-1}.$}
\label{table1}
\end{table}

\begin{table}
\begin{center}
\begin{tabular}{|c|c|c|}\hline
Method  & estimator standard deviations \\ \hline
$\varepsilon_N =N^{-\frac{5}{4}}$  & $2^{-11.25},2^{-12.50},2^{-13.75},2^{-15.00},2^{-16.25}$\\ \hline
MC and Diff. approx  & $2^{-11.27}, 2^{-12.51}, 2^{-13.75}, 2^{-15.00}, 2^{-16.25} $\\ \hline
MC and Tau-leaping  & $ 2^{-11.26}, 2^{-12.52}, 2^{-13.76}, 2^{-15.00}, 2^{-16.25}$\\ \hline
MC and Midpoint Tau-leaping  & $ 2^{-11.26}, 2^{-12.52}, 2^{-13.76}, 2^{-15.00}, 2^{-16.25}$\\ \hline
Multilevel Diff. approx  & $2^{-11.46}, 2^{-12.63}, 2^{-13.85}, 2^{-15.06}, 2^{-16.29}$\\ \hline
Biased Multilevel Tau-leaping & $2^{-11.62}, 2^{-12.81}, 2^{-13.99}, 2^{-15.19}, 2^{-16.41}$\\ \hline
Unbiased Multilevel Tau-leaping & $ 2^{-11.34}, 2^{-12.57}, 2^{-13.79}, 2^{-15.03}, 2^{-16.26}$\\ \hline
\end{tabular}
\end{center}
\caption{Actual estimator standard deviations when $\varepsilon_N = N^{-5/4}.$}
\label{table2}
\end{table}

The specifics of the implementations and results for the different Monte Carlo methods are detailed below.
\vspace{.1in}

\noindent \textbf{Diffusion Approximation plus Monte Carlo.} We took a time step of size $h = \varepsilon_N$ to generate our independent samples.  See Figure \ref{fig:example_simple}, where the best fit line is $y=1.94x  -0.88$, and Figure \ref{fig:example_54}, where the best fit line is $y=2.73x -1.37$, which  are consistent with the exponent $\alpha$ in Table~\ref{tab1}.

\vspace{.1in}

\noindent \textbf{Monte Carlo Tau-Leaping.} We took a time step of size $h = \varepsilon_N$ to generate our independent samples.  See Figure \ref{fig:example_simple}, where the best fit line is $y = 1.96x  -1.02$, and Figure \ref{fig:example_54}, where the best fit line is $y=2.76x -1.63$, which  are consistent with the exponent $\alpha$ in Table~\ref{tab1}.

\vspace{.1in}

\noindent \textbf{Monte Carlo Midpoint Tau-Leaping.}  We took a time step of size $h = \sqrt{\varepsilon_N}$.  See Figure \ref{fig:example_simple}, where the best fit line is $y = 1.44  -0.86$, and Figure \ref{fig:example_54}, where the best fit line is $y=2.10x -3.53$, which  are consistent with the exponent $\alpha$ in Table~\ref{tab1}.

\vspace{.1in}

Our implementation of the multilevel methods proceeded as follows.  We chose $h_\ell = 2^{-\ell}$ and for  $\varepsilon_N>0$ we fixed $h_L  = \varepsilon_N$ and $L = \lceil\log(h_{L})/\log(2)\rceil$ for the biased methods.  
For each level we generated $N_0$ independent sample trajectories in order to estimate $\delta_{N, \ell}$, as defined in section \ref{sec:approximations}.  Then we selected
\[
n_{\ell}=\Bigg\lceil\varepsilon_N^{-2} \sqrt{ \delta_{N,\ell} h_{\ell}}\sum_{j = 0}^L\sqrt{ \frac{\delta_{N,j}}{h_j}}\Bigg\rceil + 1, \qquad \text{ for } \ell \in\{0,1,2,\dots,L\} ,
\]
to ensure the overall variance is below the target $\varepsilon_N^2$.

\vspace{.1in}
\noindent \textbf{Multi-Level Monte Carlo Diffusion Approximation}  We used $N_0= 400$ for our pre-calculation of the variances.   See Figure \ref{fig:example_simple}, where the best fit line is $y =  0.99x  +2.75$, and Figure \ref{fig:example_54}, where the best fit line is $y=1.45x + 2.61$, which  are consistent with the exponent $\alpha$ in Table~\ref{tab1}.

\vspace{.1in}
\noindent \textbf{Multi-Level Monte Carlo Tau-Leaping.}  We used $N_0= 100$ for our pre-calculation of the variances. See Figure \ref{fig:example_simple}, where the best fit line is $y =  1.12x  +3.70$, and Figure \ref{fig:example_54}, where the best fit line is $y=1.56x + 4.64$, which  are, up to a log factor, consistent with the exponent $\alpha$ in Table~\ref{tab1}.

\vspace{.1in}

\noindent \textbf{Unbiased Tau-leaping multilevel Monte Carlo.}  For our implementation of  unbiased multilevel  tau-leaping, we set $h_L = \frac{2}{N} \text{LambertW}\left(\frac{N}{2}\right)$ and $L = \lceil\log(h_{L})/\log(2)\rceil$.  
For each level we utilized $N_0= 100$ independent sample trajectories in order to estimate $\delta_{N, \ell}, C_{\ell}, \delta_{N, E},$ and $C_E$, as defined in section \ref{sec:approximations}.  We then selected  
\[
n_{\ell}=\Bigg\lceil\varepsilon_N^{-2} \sqrt{ \frac{\delta_{N,\ell}}{C_{\ell}}}\left(\sum_{\ell = 0}^L\sqrt{ \delta_{N,\ell}C_{\ell}} +\sqrt{ \delta_{N,E}C_{E}}\right)\Bigg\rceil + 1, \qquad \text{ for } \ell \in\{0,1,2,\dots,L\},
\]
and 
\[
n_{E}=\Bigg\lceil\varepsilon_N^{-2} \sqrt{ \frac{\delta_{N,E}}{C_{E}}}\left(\sum_{\ell = 0}^L\sqrt{ \delta_{N,\ell}C_{\ell}} +\sqrt{ \delta_{N,E}C_{E}} \right)\Bigg\rceil + 1,
\]
to ensure the overall estimator variance is below our target $\varepsilon_N^2$.
See Figure \ref{fig:example_simple}, where the best fit line is $y = 1.08x  +3.71$,  and Figure \ref{fig:example_54}, where the best fit line is $y = 1.68x  +2.65$, which  are, up to a log factor, consistent with the exponent $\alpha$ in Table~\ref{tab1}.

We also used the unbiased tau-leaping multilevel Monte Carlo method with $h_L = N^{-1}$ to estimate $\E[X_1(1)]$ to accuracy $\varepsilon_N = N^{-\alpha}$, for both $\alpha=1$ and $\alpha=5/4$.  See Figures \ref{fig:example_compare1} and \ref{fig:example_compare54} for log-log plots of the required complexity when $h_L = N^{-1}$ and $h_L =\frac{2}{(\log 2)^2 N} \text{LambertW}\left(\frac{N}{2/(\log 2)^2}\right) $. As predicted in section \ref{subsec:tauMLMC}, the complexity required when $h_L = \frac{2}{(\log 2)^2 N} \text{LambertW}\left(\frac{N}{2/(\log 2)^2}\right) $ is lower by some constant factor.

 \begin{figure}
\centering
     \includegraphics[width = \textwidth]{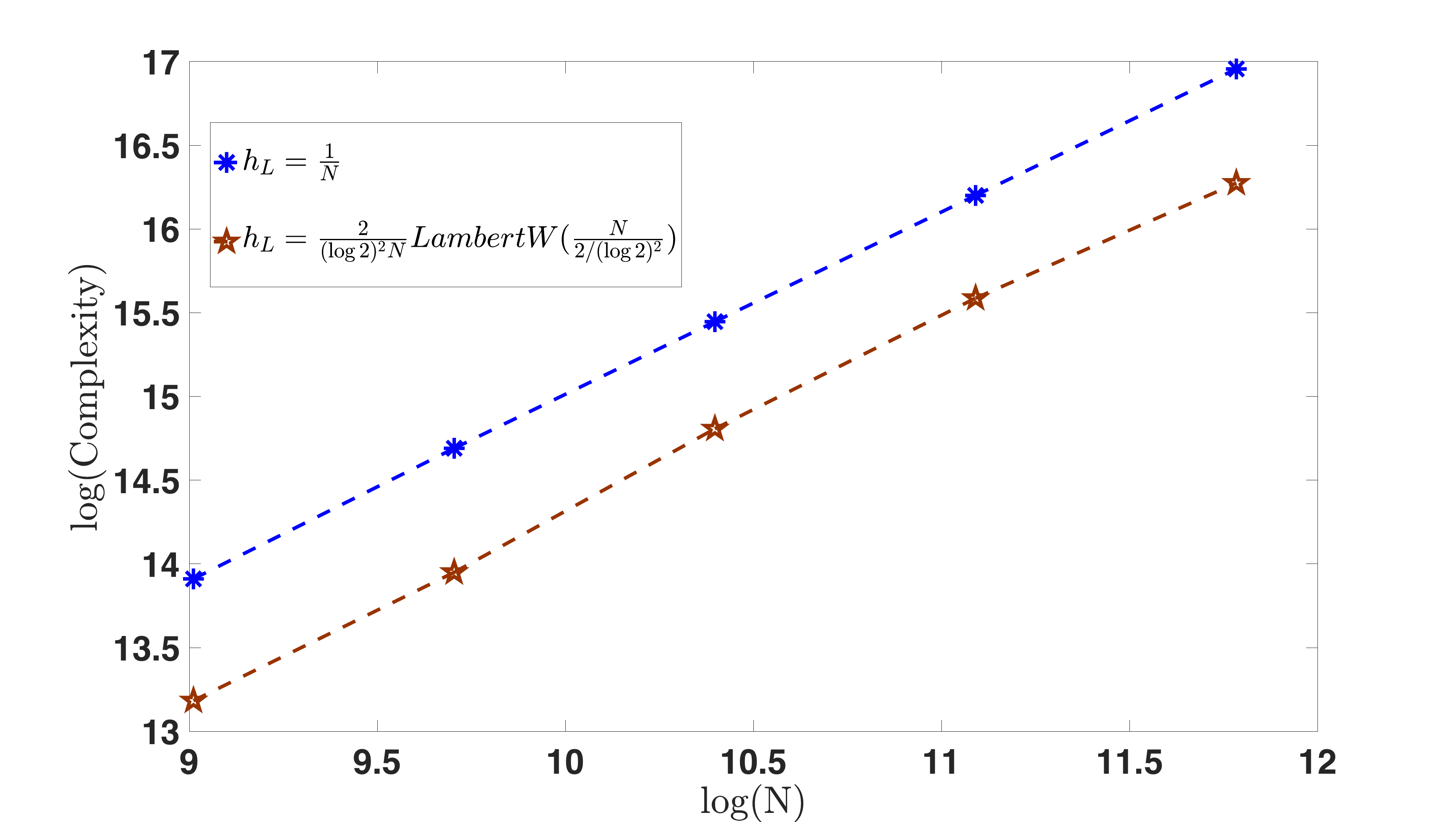}
	\caption{Complexity comparison of unbiased multilevel Monte Carlo tau-leaping when $h_L = \frac{1}{N}$  and  $h_L =\frac{2}{(\log 2)^2 N} \text{LambertW}\left(\frac{N}{2/(\log 2)^2}\right) $, with $\varepsilon_N = N^{-1}$. }
	\label{fig:example_compare1}
\end{figure}

 \begin{figure}
\centering
     \includegraphics[width = \textwidth]{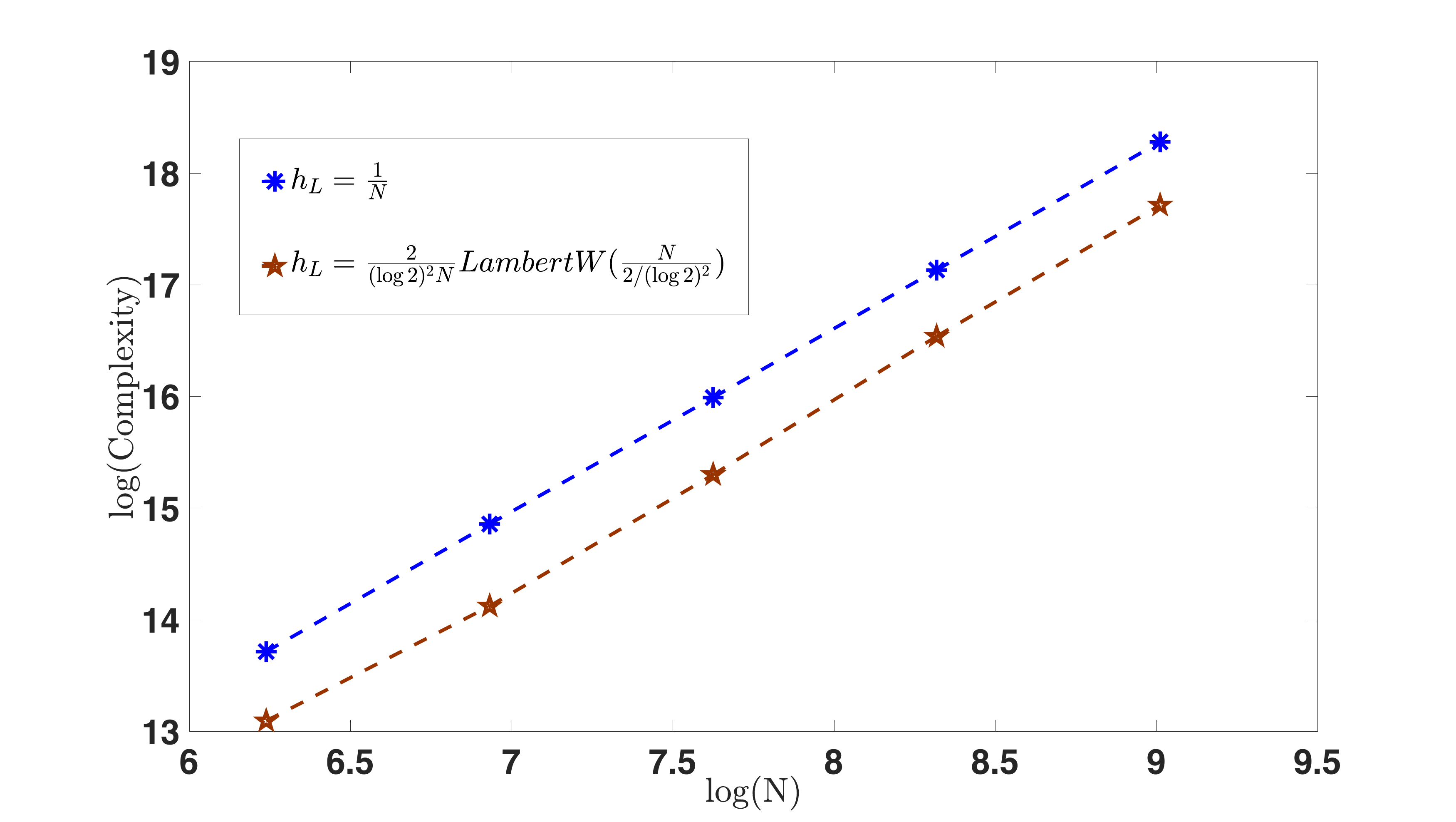}
	\caption{Complexity comparison of unbiased multilevel Monte Carlo tau-leaping when $h_L = \frac{1}{N}$  and  $h_L =\frac{2}{(\log 2)^2 N} \text{LambertW}\left(\frac{N}{2/(\log 2)^2}\right) $, with $\varepsilon_N = N^{-\frac{5}{4}}$.}
	\label{fig:example_compare54}
\end{figure}

\label{example:simple}
\end{example}

 \section{Conclusions}
\label{sec:conc}

Many researchers have 
observed in practice that approximation methods can lead to 
computational efficiency, relative to exact path simulation.
However, \emph{meaningful, rigorous justification} for  
\emph{whether} and \emph{under what circumstances}  approximation methods offer computational benefit has proved elusive.
Focusing on the classical 
scaling, we note that a useful analysis must resolve two issues:
\begin{description}
 \item[(1)] Computational complexity is most relevant for \lq\lq large\rq\rq\ problems, where 
many events take place.  However, as the system size grows the problem converges to a simpler, deterministic limit that is cheap to solve.
\item[(2)] On a fixed problem, in the traditional numerical analysis setting where mesh size tends to zero, discretization methods become arbitrarily more expensive than exact simulation  because
the exact solution is piecewise constant.
\end{description}
 
In this work, we offer what we believe to be the first rigorous complexity
analysis that allows for systematic comparison of simulation methods. 
The results, summarized in Table~\ref{tab1}, apply under the classical scaling for a family of problems
parametrized by the system size, $N$,
with accuracy requirement $N^{-\alpha}$. In this regime, we can study   
 performance 
on \lq\lq large\rq\rq\ problems when fluctuations are still relevant.

A simple conclusion from our analysis is that standard tau-leaping does offer a concrete advantage 
 over exact simulation  when the accuracy requirement is not too high, $\alpha < 1$; see the first two rows of Table~\ref{tab1}.
Also, \lq\lq second order\rq\rq\  midpoint or trapezoidal 
 tau-leaping improves on exact simulation for
$\alpha < 2$; 
row three of Table~\ref{tab1}.
  Furthermore, in this framework, we were able to analyze the use of a diffusion, or Langevin,
 approximation and the multilevel Monte Carlo versions of tau-leaping and diffusion simulation.
 Our overall conclusion is that in this scaling regime, using exact samples alone is never worthwhile.  For low accuracy 
 ($\alpha < 2/3$), second order tau-leaping with standard Monte Carlo is the most 
efficient of the methods considered.  
At higher accuracy requirements,  
$\alpha > 2/3$, multilevel Monte Carlo with a diffusion approximation is best so long as the bias inherent in perturbing the model is provably lower than the desired error tolerance.  When no such analytic bounds can be achieved, multilevel versions of tau-leaping are the methods of choice.  Moreover, for high accuracy ($\alpha > 1$), the \textit{unbiased} version is the most efficient as it does not need to take a time step smaller than $\varepsilon_N$ as the biased version must.

Possibilities for further research along the lines opened up by this work include:
\begin{itemize}
\item analyzing  other methods within this framework, for example, (a) multilevel Monte Carlo
             for the diffusion approximation using discretization methods customized for 
                       small noise systems, or (b) methods that tackle the Chemical Master Equation
                       directly using large scale deterministic ODE technology \cite{Ja13,KKNS14},
 \item development of tau-leaping methods with weak order greater than two,
\item coupling  the required accuracy to the system size in other scaling regimes, for example, 
                    to study specific problem classes with multiscale structure \cite{Ball06},
\item  determining conditions on the system for when  the diffusion approximation and Euler-Maruyama scheme achieve the $\Theta(h)$ bias given in Assumption \ref{assump:bias},
\item  determining wider classes of models and functionals $f$ for which Assumptions \ref{assump:costperpath}, \ref{assump:bias}, and \ref{assump:variance} hold.  In particular, most of the results in the literature require $\lambda_k$ to be Lipschitz and for $f$ to be a scalar valued function with domain $\Z^d$ and bounded second derivatives.

\end{itemize}

 \bibliographystyle{siam}

\bibliography{MLMC_TAU}

\end{document}